\documentclass{svproc}
\usepackage{url}

\usepackage[utf8]{inputenc}
\usepackage{amsmath,amssymb,graphicx,subfigure,float,hyperref}

\begin{document}
\mainmatter              % start of a contribution
\title{The properties of the solution for Cauchy problem of a double nonlinear time-dependent parabolic equation in non-divergence form with a source or absorption}
\titlerunning{Time-dependent parabolic equation in non-divergence form}  % abbreviated title (for running head)
\author{Mersaid Aripov \and Makhmud Bobokandov}
\authorrunning{Aripov M. and Bobokandov M.} % abbreviated author list (for running head)
\tocauthor{Mersaid Aripov and Makhmud Bobokandov}
\institute{National University of Uzbekistan, 100174, Uzbekistan\\
\email{m.boboqandov@nuu.uz}\\ 
}

\maketitle              % typeset the title of the contribution

\begin{abstract}
This paper studies the properties of solutions for a double nonlinear time-dependent parabolic equation with variable density, not in divergence form with a source or absorption. The problem is formulated as a partial differential equation with a nonlinear term that depends on both the solution and time. The main results are the existence of weak solutions in suitable function spaces; regularity and positivity of solutions; asymptotic behavior of solutions as time goes to infinity; comparison principles and maximum principles for solutions. The proofs are based on comparison methods and asymptotic techniques. Some examples and applications are also given to illustrate the features of the problem.
% We would like to encourage you to list your keywords within
% the abstract section using the \keywords{...} command.
\end{abstract}

\section{Introduction}

In the domain $Q=\left\{(t,x):t \geq t_0>0, x\in R^N\right\}$, we consider the Cauchy problem to the double nonlinear time-dependent parabolic equation with variable density, not in divergence form with a source or absorption
 \begin{eqnarray}
& \left|x\right|^{-n} \partial _{t} u=u^{q} div\left(\left|x\right|^{n_{1} } u^{m-1} \left|\nabla u^{k} \right|^{p-2} \nabla u\right)+\varepsilon \left|x\right|^{-n} t^{l} u^{\beta } ,\; \left(x,t\right)\in Q
\label{t1} \\
& u\left(0,x\right)=u_{0} \left(x\right)\ge 0,\, \, \qquad x\in R^{N}
\label{t2}
\end{eqnarray}

\noindent where $k,m\ge 1,p\ge 2,0<q<1,\, \varepsilon =\pm 1$ and nonnegative $n,l,\beta $ are given numerical parameters.

\noindent The ~(\ref{t1})-~(\ref{t2}) arises in different applications \cite{bi5}-\cite{bi7}. The equation ~(\ref{t1}) might be degenerate at the points where $u=0$ and $\nabla u=0$. Therefore, in this case, we need to consider a weak solution from having a physical sense class.

\noindent The qualitative properties of the problem ~(\ref{t1})-~(\ref{t3}) depending on $\varepsilon =1$ and $\varepsilon =-1$ are different \cite{bi1}

\noindent The ~(\ref{t1}) for the particular value of numerical parameters when $\varepsilon =1$ and $\varepsilon =-1$, intensively studied by many authors, (see \cite{bi2}-\cite{bi6} and literature therein).

\noindent Investigating qualitative properties of the problem such as Fujita type global solvability, asymptotic solution, localization of solution, finite speed propagation of distribution, blow-up solution, and so on by many authors based on self-similar solutions (for example, see \cite{bi5}- \cite{bi7},\cite{bi20} and literature therein).

\noindent Martynenko and Tedeev \cite{bi5}, studied the Cauchy problems in the case $q=n=0,k=1,l=0,\varepsilon =1$. They showed that under some restrictions on the parameters, any nontrivial solution to the Cauchy problem blows up in a finite time and established a sharp universal estimate of the solution near the blow-up point.

\noindent D. Andreucci and A.F. Tedeev studied \cite{bi18} the Neumann problem in  the case $n=l=0,k=m=1,\gamma =1$. They established a sharp universal estimate of the solution near the blow-up point and showed the condition for the solutions to exist and to which class they belonged.

\noindent Also, D. Andreucci and A.F. Tedeev proved \cite{bi19} prior supremum bounds for solutions in the case $n=l=0,\, \gamma =1,$ as $t$ approaches the time when $u$ becomes unbounded. Such bounds are universal in the sense that they do not depend on $u$.

\noindent Furthermore, recall some well-known results. In particular, when $\gamma =1,\, l=0$ authors of the work \cite{bi2}, studied the heat conduction equation with nonlinear source term. They showed that for the Cauchy problem, the critical Fujita exponent is
\[\beta _{c} =m+k(p-2)+q+\cfrac{p(1-q)}{N+n_{1} } .\] 
In the case, $l=n=0,\, k=1$ and without a source term authors \cite{bi3} proved that, depending on values of the numerical parameters and the initial value, the global solutions of the Cauchy problem exist.

\noindent The support of solutions of the equation ~(\ref{t1}) will never expand at $q>1$, while it is known that the equations in divergence form have the property of the finite (or infinite) speed of propagation of disturbance [8-10].

The aim of this paper is to study of influence a variable density and time-dependent term to an evaluation of the nonlinear processes. It is proved Fujita type global solvability, asymptotic properties of the self-similar solution based on an algorithm of nonlinear splitting. It is solved the problem with an initial approximation for the numerical solution of the problem ~(\ref{t1}), it is suggested the numerical scheme method and algorithm of solution keep nonlinear properties solution of fast diffusion, slow diffusion, critical and singular cases.

To find some self-similar solutions, which can be constructed in two ways: forward and backward, then to prove that all solutions satisfying equation ~(\ref{t1}) have the following asymptotic:
 \begin{eqnarray}\label{t3}
 f\left(\xi \right)\to C\left(a-b\xi ^{\gamma _{1} } \right)^{\gamma _{2} } \text{ at } \xi \to \xi _{b} 
\end{eqnarray}
\noindent where $\gamma _{1} =\cfrac{p}{p-1} ,\, \gamma _{2} =\cfrac{\left(p-1\right)\left(1-q\right)}{m+k(p-2)+q-1} ,\, \xi _{b} =\left\{\begin{array}{l} {\left(a/b\right)^{1/\gamma _{1} } ,\, b>0} \\ {\infty ,\, b<0} \end{array}\right. $ and $C$ is an arbitrary constant.

 We introduce the notation $v=u^{1-q} $ and put this into to problem ~(\ref{t1})-~(\ref{t2})
 \begin{eqnarray}
& r^{-n} v_{t} =r^{1-N} \cfrac{\partial }{\partial r} \bigg(r^{n_{1} +N-1} v^{m_{2} -1} \bigg|\cfrac{\partial v^{k_{2} } }{\partial r} \bigg|^{p-2} \cfrac{\partial v}{\partial r} \bigg)+\varepsilon (1-q)r^{-n} t^{l} v^{\beta _{2} } 
\label{t4} \\
& v\left|{}_{t=0} =v_{0} \left(x\right)\right. 
\label{t5}
\end{eqnarray}
\noindent where $r=|x|,m_{2} =\cfrac{m}{1-q} ,k_{2} =\cfrac{k}{1-q} ,\beta _{2} =\cfrac{\beta -q}{1-q} $.

\noindent We are looking for the solution $v(t,r)$, that has the following form 
\[v(t,r)=\bar{v}(t)w(\tau (t),\varphi (r))\] 
where 
\begin{math}
\bar{v}(t)=
\begin{cases}
{l_{1} t^{\frac{1+l}{1-\beta _{2} } } ;l\ne -1,\beta _{2} \ne 1,} \\ 
{l_{2} \left(\ln{t}\right)^{\frac{1}{1-\beta _{2} } } ;l=-1,\beta _{2} \ne 1,} \\ 
{e^{-l_{3} t^{l+1} } ;l\ne -1,\beta _{2} =1,} \\ 
{t^{-l_{4} } ;\beta _{2} =-l=1.}
\end{cases}
\end{math}

\noindent here $l_{4} =\varepsilon (1-q),l_{3} =\cfrac{l_{4} }{1+l} ,l_{2} =\left[l_{4} (\beta _{2} -1)\right]^{\frac{1}{1-\beta _{2} } } ,l_{1} =\left[l_{3} (\beta _{2} -1)\right]^{\frac{1}{1-\beta _{2} } } $.
\[
\tau (t)=
\begin{cases}
{\cfrac{l_{6} (1-\beta _{2} )t^{\frac{l_{5} }{1-\beta _{2} } +1} }{l_{5} +1-\beta _{2} } ;\, \, \, l_{5} \ne \beta _{2} -1,}
\\ {l_{6} \ln{t};\, \, \, l_{5} =\beta _{2} -1,}
\end{cases}
\] 
here $l_{5} =\left(1+l\right)\left(m_{2} +k_{2} \left(p-2\right)-1\right),l_{6} =l_{1}^{m_{2} +k_{2} (p-2)-1} .$
\[
\varphi (r)=
\begin{cases}
{\cfrac{pr^{\frac{p-n-n_{1} }{p} } }{p-n-n_{1} } ,\, \, \text{ if }\, p\ne n+n_{1} ,} \\
{\ln{r},\, \, \, \text{ if }\, \, \, p=n+n_{1} .}
\end{cases}
\] 
It is easy to check that for an unknown function $w$, satisfying the following equation
\begin{eqnarray}
\cfrac{\partial w}{\partial \tau } =\varphi ^{1-s} \cfrac{\partial }{\partial \varphi } \bigg(\varphi ^{s-1} w^{m_{2} -1} \bigg|\cfrac{\partial w^{k_{2} } }{\partial \varphi } \bigg|^{p-2} \cfrac{\partial w}{\partial \varphi } \bigg)+\cfrac{w^{\beta _{2} } +w}{l_{7} \tau }  
\label{t6}
\end{eqnarray}

\noindent where $s=\cfrac{p\left(N-n\right)}{p-n-n_{1} } ,l_{7} =-\cfrac{l_{5} +1-\beta _{2} }{1+l} ,\, n<N,\, p>n+n_{1} $.

 Put in ~(\ref{t6})
\begin{eqnarray}
w\left(\tau ,\varphi \right)=f\left(\xi \right),\, \xi =\varphi \tau ^{-1/p} 
\label{t7}
\end{eqnarray}
\noindent Then, it is easy to see that $f\left(\xi \right)$ satisfies to the following degenerate type self-similarly differential equation
\begin{eqnarray}
Af\equiv \xi ^{1-s} \cfrac{d}{d\xi } \bigg(\xi ^{s-1} f^{m_{2} -1} \bigg|\cfrac{df^{k_{2} } }{d\xi } \bigg|^{p-2} \cfrac{df}{d\xi } \bigg)+\cfrac{\xi }{p} \cfrac{df}{d\xi } +\cfrac{f+f^{\beta _{2} } }{l_{7} } =0
\label{t8}
\end{eqnarray}

\section{Asymptotic of compactly supported weak solution}

\noindent In accordance with the statement of the original problem we will consider nontrivial, nonnegative solutions of the equation ~(\ref{t8}) satisfying the following conditions:
\begin{eqnarray}
f'\left(0\right)=0,\, \, f\left(\infty \right)=0.
\label{t9}
\end{eqnarray}
\noindent Consider the following function
\begin{eqnarray}
\overline{f}\left(\xi \right)=\left(a-b\xi ^{\gamma _{1} } \right)_{+}^{\gamma _{2} } ,\, \, b=\left[\gamma _{1} \gamma _{2} pk_{2}^{p-2} \right]^{-\frac{1}{p-1} } 
\label{t10}
\end{eqnarray}
\noindent The solution of this form first was found by Zeldovich, Kompaneets, and Barenblatt in 1950 for porous medium equation \cite{bi12}. Hence, these types of solutions usually are named ZKB solutions.

\noindent Notation $\left(d\right)_{+} =\max \left\{d,0\right\}$ is used to show that we are searching the solution with compact support for initial value problem ~(\ref{t7})-~(\ref{t8}).

\noindent The constants $\gamma _{1} ,\, \gamma _{2} $ will be defined after substituting ~(\ref{t10}) into to ~(\ref{t9}) which yields
\begin{align*}
 & A\overline{f}\equiv \left(b\gamma _{1} \right)^{p} \gamma _{2}^{p-1} k_{2}^{p-2} \big[\gamma _{2} \left(m_{2} +k_{2} \left(p-2\right)-1\right)+1-p\big]\xi ^{s-1+\left(\gamma _{1} -1\right)p} \times 
\\ 
& \times \Big(a-b\xi ^{\gamma _{1} } \Big)_{+}^{\gamma _{2} \left(m_{2} +k_{2} \left(p-2\right)-1\right)-p} - \left(b\gamma _{1} \gamma _{2} \right)^{p-1} k_{2}^{p-2} \big[s-1+\left(\gamma _{1} -1\right)\left(p-1\right)\big] \times 
\\
& \times \xi ^{s-2+\left(\gamma _{1} -1\right)\left(p-1\right)} \Big(a-b\xi ^{\gamma _{1} } \Big)_{+}^{\gamma _{2} \left(m_{2} +k_{2} \left(p-2\right)-1\right)+1-p} -\cfrac{b\gamma _{1} \gamma _{2} }{p} \xi ^{s-1+\gamma _{1} } \Big(a-b\xi ^{\gamma _{1} } \Big)_{+}^{\gamma _{2} -1}+
\\ 
&  +\cfrac{\xi ^{s-1} \Big(a-b\xi ^{\gamma _{1} } \Big)_{+}^{\gamma _{2} } \Big(1+\Big(a-b\xi ^{\gamma _{1} } \Big)_{+}^{\gamma _{2} \beta _{2} -1} \Big)}{l_{7} } =0   
\end{align*}
Considering that the first and the third terms are as the second and the fourth ones respectively, we can define $\gamma _{1,2} $ from a solution of the system
\[
\begin{cases}
s-1+\left(\gamma _{1} -1\right)p=s-1+\gamma _{1} 
\\ 
\gamma _{2} \left(m_{2} +k_{2} \left(p-2\right)-1\right)-p=\gamma _{2} -1
\end{cases}
\] 
The solution of this system is 
\[
\gamma _{1} =\cfrac{p}{p-1} ,\, \, \gamma _{2} =\cfrac{p-1}{m_{2} +k_{2} \left(p-2\right)-1} .
\] 
Then we came to the following form
\begin{eqnarray}
A\overline{f}\equiv \xi ^{s-1} \left(a-b\xi ^{\gamma _{1} } \right)_{+}^{\gamma _{2} } \left[\cfrac{1}{l_{7} } -\cfrac{s}{p} +\cfrac{\left(a-b\xi ^{\gamma _{1} } \right)_{+}^{\gamma _{2} \beta _{2} -1} }{l_{7} } \right]
\label{t11}
\end{eqnarray}
\noindent Since, $\forall \xi ,\, \, \xi ^{s-1} \Big(a-b\xi ^{\gamma _{1} } \Big)_{+}^{\gamma _{2} } \ge 0$, 

\noindent we have
\[
\cfrac{1}{l_{7} } -\cfrac{s}{p} +\cfrac{\left(a-b\xi ^{\gamma _{1} } \right)_{+}^{\gamma _{2} \beta _{2} -1} }{l_{7} } \le 0,
\] 
and it is enough to satisfy the following inequality
\begin{eqnarray}
0\ge \cfrac{1}{l_{7} } -\cfrac{s}{p} =-\cfrac{1+l}{\left(1+l\right)\left(m_{2} +k_{2} \left(p-2\right)-1\right)+1-\beta _{2} } -\cfrac{N-n}{p-n-n_{1} .}
\label{t12}
\end{eqnarray}
\noindent From inequality ~(\ref{t12}), we get:
\begin{eqnarray}
\beta _{2} \ge \beta _{2c} =1+\left(1+l\right)\left[m_{2} +k_{2} \left(p-2\right)+\cfrac{p-n_{1} -N}{N-n} \right]
\label{t13}
\end{eqnarray}
\noindent $\beta _{2c} $ is a critical Fujita exponent.

\noindent Hence, the inequality from ~(\ref{t13}) we get $l_{7} \le 0$. Therefore, $A\overline{f}\le 0$ for $\xi >0$.

\noindent Above we considered case $m_{2} +k_{2} \left(p-2\right)-1\ne 0$. For the critical case $m_{2} +k_{2} \left(p-2\right)-1=0$, it is easy to show that the function
\begin{eqnarray}
& z\left(t,x\right) = \overline{v}\left(t\right)\overline{f}\left(\xi \right)
\label{t14}
\\
& \overline{f}\left(\xi \right) =
\begin{cases}
{\left(a-b\xi ^{\gamma _{1} } \right)_{+}^{\gamma _{2} } ,\quad \text{ if } m_{2} +k_{2} \left(p-2\right)\ne 1}
\\ 
{e^{-b\xi ^{\gamma _{1} } } ,\quad \text{ if } m_{2} +k_{2} \left(p-2\right)=1}
\end{cases}
\label{t15}
\end{eqnarray}

\noindent where $a=const>0$ and $b$ defined above constant is super solution.

\noindent More exactly we have

\begin{theorem}\label{A1}
 Let us $m_{2} +k_{2} \left(p-2\right)-1\ge 0,\, p>n+n_{1} ,\, \, u\left(0,x\right)\le z\left(0,x\right),\, \, x\in R^{N} $. Then for solution of the problem ~(\ref{t1})-~(\ref{t2}) an estimate
 $$u\left(t,x\right)\le z\left(t,x\right) \quad \text{ in } Q $$
\noindent hold.
\end{theorem}

\begin{proof}
Theorem \ref{A1} is proved by the comparing solution method \cite{bi1}, for comparing solution is taken the function $z\left(t,x\right)$ is considered. Substituting ~(\ref{t14}) into ~(\ref{t4}) the following inequality can be obtained:
\begin{eqnarray}
\xi ^{1-s} \cfrac{d}{d\xi } \bigg(\xi ^{s-1} \overline{f}^{m_{2} -1} \bigg|\cfrac{d\overline{f}^{k_{2} } }{d\xi } \bigg|^{p-2} \cfrac{d\overline{f}}{d\xi } \bigg)+\cfrac{\xi }{p} \cfrac{d\overline{f}}{d\xi } +\cfrac{\overline{f}+\overline{f}^{\beta _{2} } }{l_{7} } \le 0  \label{t16}
\end{eqnarray}
\noindent for the function $\bar{f}(\xi )$, the inequality ~(\ref{t16}) can be rewritten as follows:
\[
\Big(a-b\xi ^{\gamma _{1} } \Big)_{+}^{\gamma _{2} } \left[\cfrac{1}{l_{7} } -\cfrac{s}{p} +\cfrac{\Big(a-b\xi ^{\gamma _{1} } \Big)_{+}^{\gamma _{2} \beta _{2} -1} }{l_{7} } \right]\le 0
\] 
or

 $$\cfrac{1}{l_{7} } -\cfrac{s}{p} +\cfrac{\left(a-b\xi ^{\gamma _{1} } \right)_{+}^{\gamma _{2} \beta _{2} -1} }{l_{7} } \le 0  \text{ in }  Q. $$

\noindent Since $a^{\gamma _{2} \left(\beta _{2} -1\right)} \ge \left(a-b\xi ^{\gamma _{1} } \right)_{+}^{\gamma _{2} \left(\beta _{2} -1\right)} $, we have $\cfrac{s}{p} \ge \cfrac{a^{\gamma _{2} \beta _{2} -1} +1}{l_{7} } $.

\noindent Then, according to the hypothesis of Theorem \ref{A1} and comparison principle, we have:

\noindent $v\left(t,x\right)\le z\left(t,x\right)$ in  $Q$, if $v_{0} (x)\le z \left(0,x\right),{\kern 1pt} {\kern 1pt} x\in R^{N} $.

\noindent The proof of Theorem \ref{A1} is completed.
\end{proof}

\section{Global solvability}

\noindent Since, the problem ~(\ref{t1})-~(\ref{t2}) is equivalent to problem ~(\ref{t4})-~(\ref{t5}), it is sufficient to solve problem ~(\ref{t4})-~(\ref{t5}). The properties of a global solvability for weak solutions of the equation ~(\ref{t4}) were proved a comparison principle \cite{bi6}. For this goal, a new equation was constructed using the standard equation method as in \cite{bi5}:
\begin{eqnarray}
v_{+} \left(t,x\right)=\bar{v}\left(t\right)\bar{f}\left(\xi \right)
\label{t17}
\end{eqnarray}
\noindent where $\bar{v}\left(t\right)$ and $\bar{f}\left(\xi \right)$ defined above.
\begin{theorem}\label{G1}
Assume that 
$$
\gamma _{2} >0,\cfrac{s}{p} \ge \cfrac{a^{\gamma _{2} \beta _{2} -1} +1}{l_{7} } ,{\kern 1pt} {\kern 1pt} v_{0} (x)\le v_{+} \left(0,x\right),{\kern 1pt} x\in R^{N} .$$
Then, for sufficiently small $v_{0} \left(x\right),$the following holds:
\begin{eqnarray}
v\left(t,x\right)\le v_{+} \left(t,x\right) \quad \text{ in } D. 
\label{t18}
\end{eqnarray}
\end{theorem}
\noindent where the function $v_{+} \left(t,x\right)$ defined above and 
\begin{math}
D=\bigg\{\left(t,x\right):t>0, \; \left|x\right|\le \bigg(\bigg(\cfrac{a}{b} \bigg)^{p-1} \bigg(\cfrac{p-n-n_{1} }{p} \bigg)^{p} \tau \bigg)^{\frac{1}{p-n-n_{1} } } \bigg\}
.
\end{math}
\begin{proof}
Theorem \ref{G1} is proved by the comparing solution method \cite{bi1}, for comparing solution is taken the function $v_{+} \left(t,x\right)$ considered. Substituting ~(\ref{t17}) into ~(\ref{t4}) the following inequality can be obtained:
\begin{eqnarray}
\xi ^{1-s} \cfrac{d}{d\xi } \bigg(\xi ^{s-1} \overline{f}^{m_{2} -1} \bigg|\cfrac{d\overline{f}^{k_{2} } }{d\xi } \bigg|^{p-2} \cfrac{d\overline{f}}{d\xi } \bigg)+\cfrac{\xi }{p} \cfrac{d\overline{f}}{d\xi } +\cfrac{\overline{f}+\overline{f}^{\beta _{2} } }{l_{7} } \le 0  \label{t19}
\end{eqnarray}

\noindent for the function $\bar{f}(\xi )$, the inequality ~(\ref{t19}) can be rewritten as follows:
\[\xi ^{s-1} \Big(a-b\xi ^{\gamma _{1} } \Big)_{+}^{\gamma _{2} } \Bigg[\cfrac{1}{l_{7} } -\cfrac{s}{p} +\cfrac{\Big(a-b\xi ^{\gamma _{1} } \Big)_{+}^{\gamma _{2} \beta _{2} -1} }{l_{7} } \Bigg]\le 0\] 
or
\\
$
\cfrac{1}{l_{7} } -\cfrac{s}{p} +\cfrac{\Big(a-b\xi ^{\gamma _{1} } \Big)_{+}^{\gamma _{2} \beta _{2} -1} }{l_{7} } \le 0 \quad \text{ in } D .
$

\noindent Since $\Big(a-b\xi ^{\gamma _{1} } \Big)_{+}^{\gamma _{2} \left(\beta _{2} -1\right)} \le a^{\gamma _{2} \left(\beta _{2} -1\right)} $, we have $\cfrac{a^{\gamma _{2} \beta _{2} -1} +1}{l_{7} } \le \cfrac{s}{p} $.

\noindent Then, according to the hypothesis of Theorem \ref{G1} and comparison principle, we have:
\[
v\left(t,x\right)\le v_{+} \left(t,x\right) \text{ in \; }  Q, \quad \text{ if \;} v_{0} (x)\le v_{+} \left(0,x\right),\quad x\in R^{N} \]

\noindent The proof of Theorem \ref{G1} is completed.
\end{proof}

\section{Asymptotic of solution in the absorption case}

\noindent In the work \cite{bi3} authors established a large time an asymptotic of solution of the problem ~(\ref{t1})-~(\ref{t2}) for the case $q=l=0,\, \varepsilon =-1,\, \beta =\beta _{c} $. They established the following asymptotic of solution for $t\sim \infty $
\[
u(t,x)\sim \left(t\ln{t}\right)^{-\frac{1}{\beta _{c} -1} } \exp \left(-\left|x\right|^{2} /t\right),\, \, \beta _{c} =1=2/N.
\] 
Generalization large time asymptotic and behavior of the front of this result for to a degenerate nonlinear parabolic equation considered in the works [3-7] Authors of these works showed the following large time asymptotic of solution of the problem ~(\ref{t1})-~(\ref{t2})
\[
u(t,x)\sim \left(\left(T+t\right)\ln{\left(T+t\right)}\right)^{-\frac{1}{\beta _{c} -1} } \exp \left(-\left|x\right|^{2} /t\right), 
\] 
to the solution of the problem for the following critical exponent case
\[
\beta _{c} =1+\left(1+l\right)\left[m_{2} +k_{2} \left(p-2\right)-1+\cfrac{p-n_{1} -n}{N-n} \right],\, n<N,p-n_{1} -n>0
\] 
and behavior a free boundary for a particular value of the numerical parameters $\left(l=1,{\rm \; }q=0\right)$.

\noindent 
\section{Slow diffusion case}

\noindent The solution ~(\ref{t15}) is also a sub solution as it satisfies the condition
\[
\xi ^{1-s} \cfrac{d}{d\xi } \bigg(\xi ^{s-1} f^{m_{2} -1} \bigg|\cfrac{df^{k_{2} } }{d\xi } \bigg|^{p-2} \cfrac{df}{d\xi } \bigg)+\cfrac{\xi }{p} \cfrac{df}{d\xi } +\cfrac{f+f^{\beta _{2} } }{l_{7} } =0
\] 
We will show that the function ~(\ref{t15}) is an asymptotic of solutions of the problem ~(\ref{t1})-~(\ref{t2}).

\noindent Now we study asymptotic behavior of the problem ~(\ref{t1})-~(\ref{t2}), $\xi \to \, \Big(\cfrac{a}{b} \Big)_{-}^{\frac{1}{\gamma _{1} } } $.

 We will seek a solution to the problem ~(\ref{t8}) in the following form
\begin{eqnarray}
f\left(\xi \right)=\overline{f}\left(\xi \right)w\left(\eta \right),\qquad \, \, \eta =-\ln \left(a-b\xi ^{\gamma _{1} } \right)
\label{t20}
\end{eqnarray}
\noindent where $w\left(\eta \right)$ unknown function. 

\noindent Since, $\eta \to +\infty $ at $\xi \to \Big(\cfrac{a}{b} \Big)_{-}^{\frac{1}{\gamma _{1} } } $.

\noindent After the transformation ~(\ref{t20}) the equation ~(\ref{t8}) becomes
\begin{eqnarray}
\left(w^{\mu } \left|Lw\right|^{p-2} Lw\right)'+a_{1} w^{\mu } \left|Lw\right|^{p-2} Lw+a_{2} Lw+a_{3} w+a_{4} w^{\beta _{2} } =0
\label{t21}
\end{eqnarray}
\noindent where
\begin{align*}
& Lw=w'-\gamma _{2} w,\, \mu =1-\left(p-1\right)\bigg(1-\cfrac{1}{\gamma _{2} } \bigg),\, a_{1} =\cfrac{s}{\gamma _{1} } a_{0} -\gamma _{2} ,
\\ 
& a_{0} =\cfrac{e^{-\eta } }{a-e^{-\eta } } ,\, a_{2} =\cfrac{\gamma _{1} \gamma _{3} }{p} ,\, \, \gamma _{3} \, =b^{-\frac{p}{\gamma _{1} } } \gamma _{1}^{-p} k_{2}^{p-2} ,\, a_{3} =\cfrac{\gamma _{3} }{l_{7} } a_{0} ,\, a_{4} =a_{3} e^{\gamma _{2} \left(1-\beta _{2} \right)\eta }.   
\end{align*}
At first, we will demonstrate that the solution of the equation ~(\ref{t8}) has a finite limit $w_{0} $ at $\eta \to +\infty $.

 Let us take the function
\[
z\left(\eta \right)=w^{\mu } \left|Lw\right|^{p-2} Lw.
\] 
The equation ~(\ref{t21}) is transformed into
\begin{eqnarray}
\begin{cases}
w'=\gamma _{2} w+w^{-\frac{\mu }{p-1} } \left|z\right|^{\gamma _{1} -2} z,
\\ 
z'=-a_{1} z-a_{2} w^{-\frac{\mu }{p-1} } \left|z\right|^{\gamma _{1} -2} z-a_{3} w-a_{4} w^{\beta _{2} }
\end{cases}
\label{t22}
\end{eqnarray}
\begin{lemma}\label{L1}
Assume that $0<K_{1} \le K_{2} ,\, \theta _{1} <\theta _{2} \lambda \le 0,$ and let $\left(w_{1} ,z_{1} \right),\, \left(w_{2} ,z_{2} \right)$ be the solutions of the system ~(\ref{t22}) with the initial value conditions $w_{i} \left(\eta _{0} \right)=K_{i} ,\, z_{i} \left(\eta _{0} \right)=\theta _{i} \, \left(i=1,2\right)$. If $w_{1} $ and $w_{2} $ are positive in $\left[\eta _{0} ,\, +\infty \right)$ then $w_{1} \left(\eta \right)\le w_{2} \left(\eta \right),\, z_{1} \left(\eta \right)<z_{2} \left(\eta \right)$ for any $\eta \in \left[\eta _{0} ,\, +\infty \right)$.
\end{lemma}
\begin{proof}
As $0<w_{1} \left(\eta _{0} \right)\le w_{2} \left(\eta _{0} \right),\, z_{1} \left(\eta _{0} \right)<z_{2} \left(\eta _{0} \right)\le 0$,
\[\begin{array}{l} {w'_{2} \left(\eta _{0} \right)=\gamma _{2} w_{2} \left(\eta _{0} \right)+w_{2}^{-\frac{\mu }{p-1} } \left|z_{2} \left(\eta _{0} \right)\right|^{\gamma _{1} -2} z_{2} \left(\eta _{0} \right)>\gamma _{2} w_{1} \left(\eta _{0} \right)+} \\ {+w_{1}^{-\frac{\mu }{p-1} } \left|z_{1} \left(\eta _{0} \right)\right|^{\gamma _{1} -2} z_{1} \left(\eta _{0} \right)=w'_{1} \left(\eta _{0} \right)\, .} \end{array}\] 
We came to $w'_{2} \left(\eta _{0} \right)>w'_{1} \left(\eta _{0} \right)$. Then there is must exist a constant $\delta >0$ such that $w_{1} \left(\eta \right)\le w_{2} \left(\eta \right),\, z_{1} \left(\eta \right)\le z_{2} \left(\eta \right)$ on $\left[\eta _{0} ,\, \eta _{0} +\delta \right]$. By repeating this process many times, we can conclude that $w_{1} \left(\eta \right)\le w_{2} \left(\eta \right),\, z_{1} \left(\eta \right)\le z_{2} \left(\eta \right)$ conditions true for all $\eta \in \left[\eta _{0} ,\, +\infty \right)$.
\end{proof}
\begin{lemma}\label{L2}
Assume that $0<K_{1} \le K_{2} ,\, 0\ge \theta _{1} \ge \theta _{2} ,$ and let $\left(w_{1} ,z_{1} \right),\, \left(w_{2} ,z_{2} \right)$ be the solutions of the system ~(\ref{t22}) with the initial value conditions $w_{i} \left(\eta _{0} \right)=K_{i} ,\, z_{i} \left(\eta _{0} \right)=\theta _{i} \, \left(i=1,2\right)$. If $w_{1} $ and $w_{2} $ are positive in $\left[\eta _{0} ,\, +\infty \right)$ then $w_{1} \left(\eta \right)\le w_{2} \left(\eta \right),\, z_{1} \left(\eta \right)\ge z_{2} \left(\eta \right)$ for any $\eta \in \left[\eta _{0} ,\, +\infty \right)$.
\end{lemma}
\begin{proof}
From the hypotheses we have
\[\begin{array}{l} {z_{2} '\left(\eta _{0} \right)+a_{1} z_{2} \left(\eta _{0} \right)=-a_{2} w_{2}^{-\frac{\mu }{p-1} } \left|z_{2} \right|^{\gamma _{1} -2} z_{2} -a_{3} w_{2} -a_{4} w_{2}^{\beta _{2} } \left|{}_{\eta =\eta _{0} } \right. =z_{2} '\left(\eta _{0} \right)} \\ {=-a_{2} \left(w_{2}^{-\frac{\mu }{p-1} } \left|z_{2} \right|^{\gamma _{1} -2} z_{2} +\gamma _{2} w_{2} \right)-\left(a_{3} -a_{2} \right)w_{2} -a_{4} w_{2}^{\beta _{2} } \left(\eta _{0} \right)=-a_{2} w_{2} '\left(\eta _{0} \right)-} \\ {-\left(a_{3} -a_{2} \right)w_{2} -a_{4} w_{2}^{\beta _{2} } \left(\eta _{0} \right)<-a_{2} w_{1}^{-\frac{\mu }{p-1} } \left|z_{1} \right|^{\gamma _{1} -2} z_{1} -a_{3} w_{1} -a_{4} w_{1}^{\beta _{2} } \left|{}_{\eta =\eta _{0} } \right. =} \\ {=-a_{2} w_{1} '\left(\eta _{0} \right)-\left(a_{3} -a_{2} \right)w_{1} -a_{4} w_{1}^{\beta _{2} } \left(\eta _{0} \right)=z_{1} '\left(\eta _{0} \right)+a_{1} z_{1} \left(\eta _{0} \right)=z_{1} '\left(\eta _{0} \right)} \end{array}\] 
Which means that $z_{2} '\left(\eta _{0} \right)<z_{1} '\left(\eta _{0} \right),\, w_{2} '\left(\eta _{0} \right)>w_{1} '\left(\eta _{0} \right)$. Then considering the proof of Lemma \ref{L1} we have $w_{1} \left(\eta \right)\le w_{2} \left(\eta \right),\, z_{1} \left(\eta \right)>z_{2} \left(\eta \right)$ for all $\eta \in \left[\eta _{0} ,\, +\infty \right)$.
\end{proof}
\begin{theorem}\label{S1}
Let $\gamma _{2} >0$, then finite solution of the problem ~(\ref{t8})-~(\ref{t9}) has an asymptotic $f\left(\xi \right)=C\overline{f}\left(\xi \right)\left(1+o\left(1\right)\right)$ at $\xi \to \bigg(\cfrac{a}{b} \bigg)_{-}^{\frac{1}{\gamma _{1} } } $.
\end{theorem}
\begin{proof}
 We first show that the solution of the system ~(\ref{t22}) has a finite limit $w_{0} $ at $\eta \to \, +\infty $. We know that the function $w_{\eta } $ is bounded. So it is enough to prove that it is monotone non-increasing in $\left[\eta _{0} ,\, +\infty \right)$. Let us take $w\equiv 1$, we can see that it is a sub-solution from ~(\ref{t20}). Then for other solutions $w_{1} \left(\eta \right)$, we have $w_{1} '\left(\eta \right)\le 0$. That means that any solution is non-increasing in $\left[\eta _{0} ,\, \eta _{1} \right)$, for any $\eta _{1} >\eta _{0} $, where difference $\eta _{1} -\eta _{0} $ sufficiently small. By considering Lemma \ref{L1} and Lemma \ref{L2} we can find two solutions $w1,{\rm \; }w2$ such that $w_{1} \left(\eta _{1} \right)=w_{2} \left(\eta _{0} \right)$. We can conclude
\[
w_{1} \left(\eta _{1} \right)\ge w_{1} \left(\eta _{0} \right).
\] 
From the arbitrariness of $\eta _{1} $, we see that $w_{1} $ is monotonic in $\left[\eta _{0} ,\, +\infty \right)$. Thus, it has limit at $\eta \to \, +\infty $. Here taking into consideration that
\[
{\mathop{\lim }\limits_{\eta \to \, +\infty }} \cfrac{be^{-\eta } }{a-e^{-\eta } } \, \to \, 0,\, w'=0,
\] 
Hence, we get the following
\[
{\mathop{\lim }\limits_{\eta \to \, +\infty }} a_{0} ={\mathop{\lim }\limits_{\eta \to \, +\infty }} a_{3} =0,\,
{\mathop{\lim }\limits_{\eta \to \, +\infty }} a_{1} =-\gamma _{2} ,\, \,
{\mathop{a}\limits^{0}} _{4} ={\mathop{\lim }\limits_{\eta \to \, +\infty }} a_{4} =
\begin{cases} 
{\cfrac{\gamma _{3} }{al_{7} } ,\, \, 
\gamma _{2} \left(1-\beta _{2} \right)=1} \\
{0,\, \, \gamma _{2} \left(1-\beta _{2} \right)<1} \\
{+\infty ,\, \, \gamma _{2} \left(1-\beta _{2} \right)>1} 
\end{cases}
\] 
In ~(\ref{t21}) we will get the following algebraic equation at $\xi \to \bigg(\cfrac{a}{b} \bigg)_{-}^{\frac{1}{\gamma _{1} } } $,
\[\gamma _{2}^{p} w_{0}^{m_{2} k_{2} \left(p-2\right)-1} -\cfrac{b^{-\frac{p}{\gamma _{1} } } \gamma _{1}^{1-p} \gamma _{2} k_{2}^{p-2} }{p} +{\mathop{a}\limits^{0}} _{4} w_{0}^{\beta _{2} -1} =0\] 
Considering ~(\ref{t20}) the theorem has been proved.
\end{proof}

\section{Fast diffusion case}

\noindent In this case we have diﬀerent families of solutions, oscillating near at $ \xi =+\infty $, see for example \cite{bi14}-\cite{bi15}. For every one of these families, there is a constant $C$ from ~(\ref{t3}). These kinds of solutions are called eigenfunction on nonlinear media \cite{bi16}. The upper solution for the problem ~(\ref{t1})-~(\ref{t2}) is obtained by nonlinear splitting \cite{bi17}.
\[
g\left(\xi \right)=A\left(a+\xi ^{\gamma _{1} } \right)^{\gamma _{2} } ,
\] 
where 
$
A=\left[\cfrac{\left(\gamma _{1} \gamma _{2} \right)^{1-p} k_{2}^{2-p} }{p} \right]^{\frac{\gamma _{2} }{p-1} } ,\, a=const>0.
$

\begin{theorem}
Let $\gamma _{2} \left(\beta _{2} -1\right)\le 0$, then finite solution of the problem ~(\ref{t8})-~(\ref{t9}) has an asymptotic at $f\left(\xi \right)=Cg\left(\xi \right)$.
\end{theorem}
\begin{proof}
Let us introduce a new transformation and use it to ~(\ref{t8})
\begin{eqnarray}
f\left(\xi \right)=g\left(\xi \right)w\left(\eta \right),\, \, \, \, \eta =\ln \left(a+\xi ^{\gamma _{1} } \right).
\label{t23}
\end{eqnarray}
\noindent Admit that
\[\eta \xrightarrow{\xi \to \, +\infty } +\infty .\] 
After the transformation ~(\ref{t23}), we have
\begin{eqnarray}
\left(w^{\mu } \left|Lw\right|^{p-2} Lw\right)'+b_{1} w^{\mu } \left|Lw\right|^{p-2} Lw+b_{2} Lw+b_{3} w+b_{4} w^{\beta _{2} } =0
\label{t24}
\end{eqnarray}
\noindent where  
\begin{eqnarray}
& Lw=w'+\gamma _{2} w,\, 
b_{1} =\cfrac{s}{\gamma _{1} } b_{0} +\gamma _{2} ,
b_{0} =\cfrac{1}{1-ae^{-\eta } } ,
\\ \notag 
& b_{2} =\cfrac{\gamma _{1}^{1-p} A^{\frac{1-p}{\gamma _{2} } } }{p} ,\, 
b_{3} =\cfrac{\gamma _{1}^{-p} A^{\frac{1-p}{\gamma _{2} } } }{l_{7} } b_{0} ,\, 
b_{4} =\cfrac{\gamma _{1}^{-p} A^{\beta _{2} -m_{2} -k_{2} \left(p-2\right)} b_{0} }{l_{7} } e^{\gamma _{2} \left(\beta _{2} -1\right)\eta }
\label{t25}
\end{eqnarray}
\noindent Taking into consideration that $w\left(\eta \right)$ has a finite limit (see the proof of theorem \ref{S1}) and
\[{\mathop{\lim }\limits_{\eta \to \, +\infty }} \cfrac{1}{1-ae^{-\eta } } \, \to \, 1,\, \, w'=0\] 
Hence, we get the following
\begin{align*}
{\, {\mathop{\lim }\limits_{\eta \to \, +\infty }} b_{1} =\cfrac{s}{\gamma _{1} } +\gamma _{2} ,\, \, 
{\mathop{\lim }\limits_{\eta \to \, +\infty }} b_{2} =\cfrac{\gamma _{1}^{1-p} A^{\frac{1-p}{\gamma _{2} } } }{p} ,\, \, 
{\mathop{\lim }\limits_{\eta \to \, +\infty }} b_{3} =\cfrac{\gamma _{1}^{-p} A^{\frac{1-p}{\gamma _{2} } } }{l_{7} } }
\\ 
{\mathop{b}\limits^{0}} _{4} ={\mathop{\lim }\limits_{\eta \to \, +\infty }} b_{4} =
\begin{cases} 
{\cfrac{\gamma _{1}^{-p} A^{\beta _{2} -m_{2} -k_{2} \left(p-2\right)} }{l_{7} } ,\, \, \gamma _{2} \left(\beta _{2} -1\right)=0}
\\ 
{0,\, \, \gamma _{2} \left(\beta _{2} -1\right)<0} \\ {+\infty ,\, \, \gamma _{2} \left(\beta _{2} -1\right)>0} 
\end{cases}
\end{align*}
We come to the following algebraic equation
\[\left(\cfrac{s}{\gamma _{1} } +\gamma _{2} \right)\gamma _{2}^{p-1} C^{\mu +p-2} +A^{\frac{1-p}{\gamma _{2} } } \gamma _{1}^{-p} \left(\cfrac{\gamma _{1} \gamma _{2} }{p} +\cfrac{1}{l_{7} } \right)+{\mathop{b}\limits^{0}} _{4} C^{\beta _{2} -1} =0\] 
at $\eta \to \, +\infty $. On behalf of ~(\ref{t25}) we have proved that $f\left(\xi \right)=Cg\left(\xi \right)\left(1+o\left(1\right)\right)$, where $C$ is the root of last algebraic equation.
\end{proof}

\section{Numerical analysis of solutions}

\noindent The main difficulty of numerical research for the problem ~(\ref{t1})-~(\ref{t2}) arises from the non-uniqueness of solutions. To choose appropriate an initial approximation we use asymptotic formulas established above

\noindent In slow diffusion case as an initial approximation we use $u_{0} \left(t,x\right)=l_{1} t^{\frac{1+l}{1-\beta _{2} } } \overline{f}\left(\xi \right)$.

\noindent Speeds of perturbation are finite for these cases and we can see how they slow down and maximal speed is at $t=0$.

\noindent For numerical analysis has been used tridiagonal matrix algorithm. The results of numerical experiments show quick convergence of the iteration process. The number of iterations did not exceed 3. Below listed evolution of diffusion processes finite speed of perturbation in the one-dimensional case.
\begin{figure}[H]
    \centering
    \caption{One-dimensional case}
    \subfigure[$q=0.8, m=k=1, p=3, n=0, \gamma=0 $]{\includegraphics[width=0.495\textwidth,height=5cm]{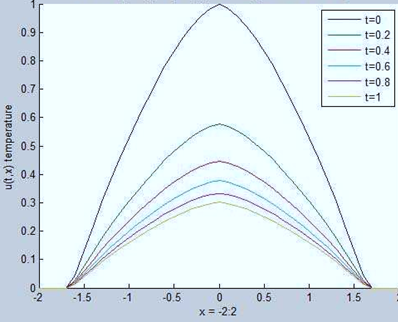}}
    \subfigure[$q=0.9, m=k=1, p=3, n=0, \gamma=0 $]{\includegraphics[width=0.495\textwidth,height=5cm]{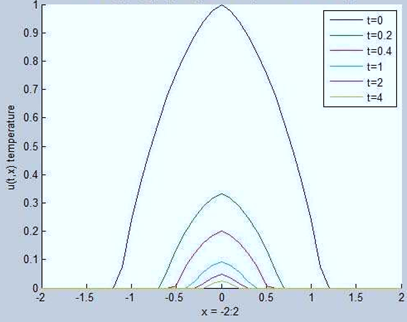}}
\end{figure}
\begin{figure}[H]
    \centering
    \caption{Two-dimensional case}
%    \\ \text{\textbf{Two-dimensional case}} \\
    \subfigure[$q=0, m=3, k=3.2, p=4, n=0.2, \gamma =1, l=0, \beta =1.3 $]{\includegraphics[width=0.495\textwidth, height=5cm]{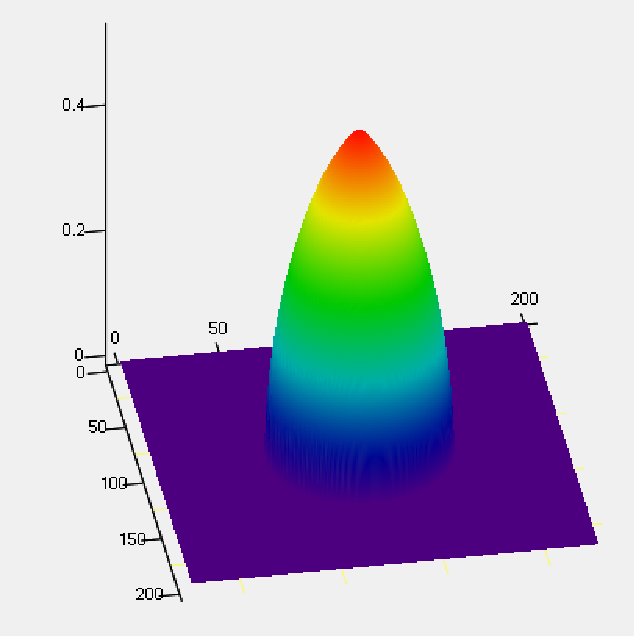}}
    \subfigure[$q=0, m=3.1, k=3, p=3.4, n=0.3, \gamma =1, l=0, \beta =3 $]{\includegraphics[width=0.495\textwidth, height=5cm]{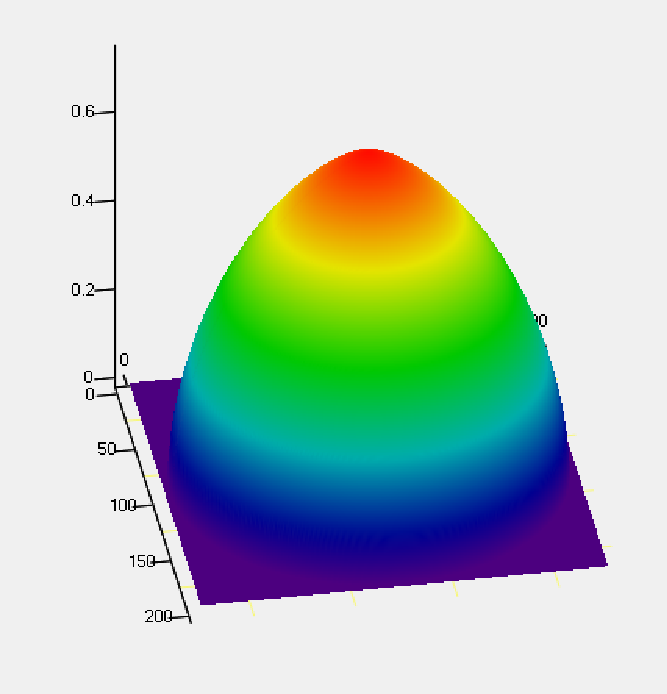}}
\end{figure}
%\clearpage

\end{document}